\documentclass[12pt,a4paper]{article}
\usepackage{amsmath,amsthm,amssymb,latexsym,epic}
\usepackage{enumerate}

\newtheorem{theorem}{Theorem}[section]
\newtheorem{lemma}[theorem]{Lemma}
\newtheorem{corollary}[theorem]{Corollary}
\newtheorem{proposition}[theorem]{Proposition}

\newtheorem{example}[theorem]{Example}
\newtheorem{question}[theorem]{Question}

\renewcommand{\phi}{\varphi}

\begin{document}

\title{On nilpotent subsemigroups in some matrix semigroups}
\author{Olexandr Ganyushkin and Volodymyr Mazorchuk}
\date{}
\maketitle

\begin{abstract}
We describe maximal nilpotent subsemigroups of a given nilpotency 
class in the semigroup $\Omega_n$ of all $n\times n$
real matrices with non-negative coefficients and the semigroup
$\mathbf{D}_n$ of all doubly stochastic real matrices.
\end{abstract}

\section{Introduction}\label{s1}

Non-negative, stochastic and doubly stochastic matrices are 
important and popular objects of study in many branches of
modern mathematics. Surprisingly enough, the algebraic
properties of the semigroup $\Omega_n$ of all $n\times n$
real matrices with non-negative coefficients, and the 
semigroups $\mathbf{P}_n$ and $\mathbf{D}_n$ of all stochastic
and doubly stochastic real matrices respectively, are
studied rather superficially. One of the possible explanations
for this fact is complexity of the algebraic structure of these
semigroups. 

Among the classical results on the algebraic structure of these
semigroups one should mention the descriptions of the maximal
subgroups in the semigroup $\mathbf{P}_n$, see \cite{Sc}, and 
the descriptions of the maximal subgroups in the semigroup 
$\mathbf{D}_n$, see \cite{Fa}. In particular, it was shown that 
every maximal subgroup in $\mathbf{P}_n$ is isomorphic to
the symmetric group $\mathrm{S}_k$ for some $k\leq n$, and that 
every maximal subgroup in $\mathbf{D}_n$ is isomorphic to
the direct product of symmetric groups. Later on there appeared
alternative and easier proofs of these results, discovered by 
different mathematicians, see for example 
\cite{Sc2,Wa1,Wa2,HM,Ta1}. In \cite{Fl} it was shown that 
each maximal subgroup of $\Omega_n$ is isomorphic to some
full monomial group of degree $k\leq n$ over the group
of positive real numbers (see also \cite{HM,Pl,Ta1,Me}).

Apart from the study of maximal subgroups, a lot of attention was
paid to the description of regular elements in these semigroups,
see \cite{Wa1,Wa2,Ta1,Pl,MP1,MP2,Py}, and Green's relations.
Green's relations were studied for regular elements in \cite{Ro}
and for arbitrary elements in \cite{Ya1,Ya2,HMP,Ta2}.
In the semigroup $\mathbf{D}_n$ Green's relations for regular
elements were studied in \cite{MP1,Wa1}, and for
arbitrary elements in \cite{MP2}. The paper \cite{Ta1}, apart from new 
results, contains several new proofs of already known facts 
about the regular elements, Green's relations and maximal subgroups
of the semigroups $\Omega_n$, $\mathbf{P}_n$ and $\mathbf{D}_n$.

Some algebraic properties of the set $E(\mathbf{D}_n)$ of idempotents
of the semigroup $\mathbf{D}_n$ were studied in \cite{SR}, and
indecomposable elements of $\Omega_n$ and $\mathbf{D}_n$ were
studied in \cite{RS,PHS}. In \cite{Sc,Ma} the minimal
ideals of $\mathbf{P}_n$ were studied. Maximal subsemigroups of
$\Omega_n$, $\mathbf{P}_n$ and $\mathbf{D}_n$ were studied in 
\cite{HM}. And, finally, in \cite{GM1,GM2} the present authors 
studied a homomorphism from the semigroup $\mathbf{D}_n$ to the 
semigroup $\mathrm{B}_n$ of all binary relations on an $n$-element set, 
the image of which has several interesting extremal properties.

In the present paper we describe maximal nilpotent subsemigroups
of a given nilpotency class in the semigroups $\Omega_n$ 
(see Section~\ref{s3}) and $\mathbf{D}_n$ (see Section~\ref{s4}). 
As a preparatory work, we develop some general reduction technique
for the study of maximal nilpotent subsemigroups via epimorphisms
in Section~\ref{s2}.
\vspace{0.2cm}

\noindent
{\bf Acknowledgments.}
The paper was written during the visit of the first author to 
Uppsala University, which was supported by The Swedish Institute.
The financial support of The Swedish Institute and the hospitality 
of Uppsala University are gratefully acknowledged. For the second 
author the research was partially supported by The Swedish 
Research Council. 

\section{Epimorphisms of nilpotent semigroups}\label{s2}

Let $S$ be a semigroup with the zero element $0$. A subsemigroup,
$T\subset S$, is called {\em nilpotent} provided that there exists
$k\in\mathbb{N}$ such that $T^k=0$. The minimal $k$ with this property is
called the {\em nilpotency class} of $T$.

For a positive integer, $m$, we denote by $\mathrm{Nil}_m(S)$ the 
set of all nilpotent subsemigroups of $S$ of nilpotency class at most 
$m$, and by $\mathrm{Nil}_m^{\mathrm{max}}(S)$~-- the set of all maximal
elements in $\mathrm{Nil}_m(S)$ with respect to inclusion. Abusing the
language, the elements of $\mathrm{Nil}_m^{\mathrm{max}}(S)$ will be called
{\em maximal nilpotent subsemigroups of nilpotency class $m$} of $S$.
Set also
\begin{displaymath}
\mathrm{Nil}(S)=\cup_{m\in\mathbb{N}}\mathrm{Nil}_m(S)
\end{displaymath}
and let $\mathrm{Nil}^{\mathrm{max}}(S)$ be the set of all maximal
elements in $\mathrm{Nil}(S)$ with respect to inclusion. The elements of
$\mathrm{Nil}^{\mathrm{max}}(S)$ are {\em maximal nilpotent 
subsemigroups} of $S$.

\begin{lemma}\label{lemma1}
Let $\varphi:S\to T$ be a surjective homomorphism of semigroups
with zero. Assume that $T$ is a nilpotent semigroup of class $n$ 
and $\varphi^{-1}(0)$ is a nilpotent semigroup of class $m$. Then
$S$ is a nilpotent semigroup of class at most $mn$.
\end{lemma}

\begin{proof}
For arbitrary $a_1,\dots,a_{nm}\in S$ we have
\begin{displaymath}
\varphi(a_1\cdots a_{n})=
\varphi(a_{n+1}\cdots a_{2n})=\dots=
\varphi(a_{n(m-1)+1}\cdots a_{nm})=0
\end{displaymath}
and hence $a_{ni+1}\cdots a_{n(i+1)}\in\varphi^{-1}(0)$
for every $i=0,1,\dots,m-1$. Therefore
\begin{displaymath}
(a_1\cdots a_{n})\cdot (a_{n+1}\cdots a_{2n})\cdots
(a_{n(m-1)+1}\cdots a_{nm})=0.
\end{displaymath}
\end{proof}

\begin{theorem}\label{theorem1}
Let $\varphi:S\to T$ be a surjective homomorphism of semigroups
with zero and assume $\mathrm{Nil}^{\mathrm{max}}(S)\neq \varnothing$. 
Then:
\begin{enumerate}[(a)]
\item\label{theorem1.1} $\varphi^{-1}$ induces a bijection
between $\mathrm{Nil}^{\mathrm{max}}(T)$ and 
$\mathrm{Nil}^{\mathrm{max}}(S)$ if and only if 
$\varphi^{-1}(0)$ is a nilpotent subsemigroup of $S$.
\item\label{theorem1.2} If $\varphi^{-1}(0)=0$, then for 
any $U\in \mathrm{Nil}(S)$ the semigroups $U$ and 
$\varphi(U)$ have the same nilpotency class; and for
any $V\in \mathrm{Nil}(T)$ the semigroups $V$ 
and $\varphi^{-1}(V)$ have
the same nilpotency class. Moreover, for every $k\in\mathbb{N}$
the maps $\varphi$ and $\varphi^{-1}$ induce mutually inverse
bijections between $\mathrm{Nil}_k^{\mathrm{max}}(S)$
and $\mathrm{Nil}_k^{\mathrm{max}}(T)$.
\end{enumerate}
\end{theorem}

\begin{proof}
We start with \eqref{theorem1.1}. The necessity is obvious. To
prove the sufficiency we observe that from Lemma~\ref{lemma1}
it follows that if $\varphi^{-1}(0)$ is nilpotent,
then  $\varphi^{-1}$ induces a map from $\mathrm{Nil}(T)$ to 
$\mathrm{Nil}(S)$, which is obviously injective. 

Let $A\in \mathrm{Nil}^{\mathrm{max}}(T)$. If
$\varphi^{-1}(A)\not\in \mathrm{Nil}^{\mathrm{max}}(S)$, then
there exists $B\in \mathrm{Nil}(S)$ such that
$\varphi^{-1}(A)\subsetneqq B$. This implies
$A\subsetneqq \varphi(B)$. On the other hand, $\varphi(B)$ is
obviously nilpotent since so is $B$. This contradicts the
maximality of $A$. Hence $\varphi^{-1}(A)\in 
\mathrm{Nil}^{\mathrm{max}}(S)$.

Let now $B\in \mathrm{Nil}^{\mathrm{max}}(S)$. Then for 
every nilpotent subsemigroup $A$ of $T$ such that $\varphi(B)\subset A$
we have $B\subset \varphi^{-1}(A)$. Since $\varphi^{-1}(A)$ is nilpotent,
from the maximality of $B$ we get $B=\varphi^{-1}(A)$ and 
$\varphi(B)=A$. This means that $\varphi(B)\in \mathrm{Nil}^{\mathrm{max}}(T)$
and the map
\begin{displaymath}
\begin{array}{ccc}
\mathrm{Nil}^{\mathrm{max}}(T) & \longrightarrow &
\mathrm{Nil}^{\mathrm{max}}(S)\\
A & \mapsto & \varphi^{-1}(A)
\end{array}
\end{displaymath}
is surjective. This proves \eqref{theorem1.1}.

Since the homomorphic image of a nilpotent semigroup is a nilpotent
semigroup, whose nilpotency class is not greater than the nilpotency
class of the original semigroup, the first part of \eqref{theorem1.2}
follows from Lemma~\ref{lemma1}. Analogously to the proof of
\eqref{theorem1.1} one shows that $\varphi$ induces the bijective map
\begin{displaymath}
\begin{array}{ccc}
\mathrm{Nil}_k^{\mathrm{max}}(T) & \longrightarrow &
\mathrm{Nil}_k^{\mathrm{max}}(S)\\
A & \mapsto & \varphi^{-1}(A),
\end{array}
\end{displaymath}
and the statement \eqref{theorem1.2} follows. This completes the proof.
\end{proof}

\section{Maximal nilpotent subsemigroups in $\Omega_n$}\label{s3}

Consider the semigroup $\mathrm{M}^0_n$ of all real $n\times n$
matrices of order $n$, each row and each column of which contain
at most one non-zero component. The semigroup $\mathrm{M}^0_n$
is naturally identified with the Rees matrix semigroup
$\mathrm{M}^0(\mathbb{R}^*;I,\Lambda;E)$, where
$I=\Lambda=\{1,2,\dots,n\}$ and the sandwich-matrix $E$ is the
identity matrix of order $n$. Since the identity matrix does not
contain any zero rows or columns, the semigroup $\mathrm{M}^0_n$
is a regular semigroup. Recall that $\Omega_n$ denotes the semigroup
of all $n\times n$ matrices with non-negative real coefficients.
Set $\tilde{\mathrm{M}}^0_n=\mathrm{M}^0_n\cap \Omega_n$. Then
$\tilde{\mathrm{M}}^0_n$ can be identified with the Rees matrix semigroup
$\mathrm{M}^0(\mathbb{R}^+;I,\Lambda;E)$, where $\mathbb{R}^+$
denotes the multiplicative group of positive real numbers. 

\begin{theorem}\label{theorem2}
Let $S$ denote one of the semigroups
$\mathrm{M}^0_n$, $\tilde{\mathrm{M}}^0_n$, or $\Omega_n$. Then
\begin{enumerate}[(a)]
\item \label{theorem2.1} The semigroup $S$ contains $n!$
maximal nilpotent subsemigroups, each of nilpotency class
$n$. These subsemigroups are in a natural bijection with 
linear orders on the set $\{1,2,\dots,n\}$.
\item \label{theorem2.2} The maximal nilpotent subsemigroup of
$S$, which corresponds to the linear order 
$i_1\prec i_2\prec\dots\prec i_n$, has the form
\begin{displaymath}
T=\{(a_{k,l})\in S:a_{k,l}\neq 0\Rightarrow k\prec l\}.
\end{displaymath}
\item \label{theorem2.3}
Let the maximal nilpotent subsemigroups $T_1$ and $T_2$ of
$S$ correspond to the linear orders 
$i_1\prec_1 i_2\prec_1\dots\prec_1 i_n$ and
$j_1\prec_2 j_2\prec_2\dots\prec_2 j_n$ respectively. Then
$T_2=M^{-1}T_1 M$, where $M=(m_{i,j})$ is the monomial matrix,
which corresponds to the permutation 
\begin{displaymath}
\pi=\left(\begin{array}{ccc}i_1&\dots&i_n\\j_1&\dots&j_n\end{array}
\right)
\end{displaymath}
(that is $M$ is a $(0,1)$-matrix such that $m_{i,j}=1$ if and only
if $j=\pi(i)$). In particular, all maximal nilpotent subsemigroups
of $S$ are isomorphic. 
\item \label{theorem2.4}
For every $k\in\mathbb{N}$, $k<n$, the semigroup $S$ contains
$\sum_{i=0}^{k-1}(-1)^i\binom{k}{i}(k-i)^n$ maximal nilpotent 
subsemigroups of nilpotency class $k$. These maximal nilpotent 
subsemigroups of nilpotency class $k$ are in a natural
bijection with decompositions of the set $\{1,2,\dots,n\}$ into
an ordered union of $k$ pairwise disjoint and non-empty blocks.
\item \label{theorem2.5}
Let $N_1\cup\dots\cup N_k=\{1,\dots,n\}$ be a decompositions 
into an ordered union of pairwise disjoint non-empty blocks.
Then the maximal nilpotent subsemigroup of nilpotency class
$k$ of $S$, which corresponds to this partition, is
\begin{displaymath}
T=\{(a_{i,j})\in S:a_{i,j}\neq 0\text{ and }i\in N_p, j\in N_q
\Rightarrow p<q\}.
\end{displaymath}
\end{enumerate}
\end{theorem}

\begin{proof}
Consider the map $\psi:(a_{i,j})\mapsto (\hat{a}_{i,j})$ from the
set $M_n(\mathbb{R})$ of all real $n\times n$ matrices to the
set of all $(0,1)$-matrices of size $n\times n$, defined via
\begin{displaymath}
\hat{a}_{i,j}=
\begin{cases}
1,& a_{i,j}\neq 0;\\
0, & \text{ otherwise. }
\end{cases}
\end{displaymath}
As $\mathrm{M}^0_n=\mathrm{M}^0(\mathbb{R}^+;I,\Lambda;E)$,
the restriction $\psi|_{\mathrm{M}^0_n}$ is a surjective
homomorphism from $\mathrm{M}^0_n$ to the multiplicative semigroup
$\mathrm{R}_n$ of those $(0,1)$-matrices of size $n\times n$, each 
row and column of which contains at most one non-zero element.
The semigroup $\mathrm{R}_n$ is usually called the {\em Rook monoid},
and it is canonically isomorphic to the symmetric inverse semigroup
$\mathcal{IS}_n$ of all partial injections on the set
$\{1,2,\dots,n\}$, see for example \cite{So}. Obviously,
$\psi^{-1}(0)=0$. Hence for the semigroup $\mathrm{M}^0_n$
the statements \eqref{theorem2.1}-\eqref{theorem2.3} follow
from Theorem~\ref{theorem1} and the description of all maximal
nilpotent subsemigroups in $\mathcal{IS}_n$, given in 
\cite{GK1}. The statements \eqref{theorem2.4}-\eqref{theorem2.5} 
follow for $\mathrm{M}^0_n$ from Theorem~\ref{theorem1} and the 
description of all maximal nilpotent subsemigroups of nilpotency
class $k$ in $\mathcal{IS}_n$, given in  \cite{GK2}. 

For the semigroup $\tilde{\mathrm{M}}^0_n$ the proof is identical.
Hence it remains to consider the case $S=\Omega_n$. If
$M\in \Omega_n$, the matrix $\psi(M)$ can be considered as 
the matrix of some binary relation on $\{1,2,\dots,n\}$. Then the
restriction $\psi|_{\Omega_n}$ defines a surjective homomorphism
from $\Omega_n$ on the semigroup $\mathrm{B}_n$ of all binary
relations on $\{1,2,\dots,n\}$. Moreover, $\psi^{-1}(0)=0$.
Therefore the statements \eqref{theorem2.1}-\eqref{theorem2.5}
for the semigroup $\Omega_n$ follow from Theorem~\ref{theorem1} 
and the description of all maximal nilpotent subsemigroups
(of a given nilpotency class) in $\mathrm{B}_n$, given 
in \cite[Theorem~5.15 and Theorem~6.1]{GM3}.
\end{proof}

\begin{corollary}\label{corollary1}
Each nilpotent matrix from $\Omega_n$ contains at most
$n(n-1)/2$ non-zero elements.
\end{corollary}

\begin{proof}
This follows from Theorem~\ref{theorem2}\eqref{theorem2.2}.
\end{proof}

From Theorem~\ref{theorem2} we obtain that the semigroup
$T$ of all upper triangular $n\times n$ matrices with non-negative
real coefficients and zero diagonal is a maximal nilpotent
subsemigroup of $\Omega_n$, and all other maximal nilpotent
subsemigroups are obtained from $T$ via conjugation with
monomial matrices. On the other hand, for every non-degenerate
matrix $A\in \Omega_n$ the semigroup $A^{-1}TA$ is
nilpotent. However, if $A$ is not monomial, then for every monomial
matrix $B$ we have $A^{-1}TA\neq B^{-1}TB$ (this follows immediately
from the obvious fact that $C^{-1}TC=T$ is possible only for
a diagonal monomial matrix $C$). Hence $A^{-1}TA$ is not a maximal
nilpotent subsemigroup of $\Omega_n$, which means that
$A^{-1}TA\not\subset\Omega_n$. By our choice of $A$ we have that
the matrix $A^{-1}MA$ can contain negative coefficients only if
there are some negative coefficients in the matrix $A^{-1}$. This
proves the following:

\begin{corollary}\label{corollary2}
The group of invertible elements of $\Omega_n$ coincides with
the complete monomial group of degree $n$ over positive reals.
\end{corollary}

We note that Corollary~\ref{corollary2} can be also derived
from the description of all maximal subgroups in $\Omega_n$,
see \cite[Theorem~1]{Fl}, \cite[Corollary~1]{Pl}, 
\cite[Corollary~3.3]{Ta1}.

\section{Maximal nilpotent subsemigroups in $\mathbf{D}_n$}\label{s4}

The structure of maximal nilpotent semigroups in the 
semigroup $\mathbf{D}_n$ of doubly stochastic matrices 
is more complicated. In \cite{GM1} it is shown that the
restriction of the map $\psi$ from the proof of Theorem~\ref{theorem2}
to $\mathbf{D}_n$ is a surjective homomorphism on the factor power
$\mathcal{FP}^+(\mathrm{S}_n)$ of the symmetric group
$\mathrm{S}_n$. However, the zero element of the semigroup
$\mathcal{FP}^+(\mathrm{S}_n)$ (the latter being considered 
as a subsemigroup of $\mathrm{B}_n$) is the full relation,
and hence $\psi^{-1}(0)$ coincides with the semigroup of all
doubly stochastic matrices, all coefficients of which are positive.
However, this subsemigroup of $\mathbf{D}_n$ is not nilpotent.
Indeed, the zero element of $\mathbf{D}_n$ is the matrix
$\mathbb{O}_n$, all coefficients of which are equal to $1/n$.
The matrix $A_t=t\mathbb{O}_n+(1-t)E$, where $0<t<1$ is doubly stochastic
and has positive coefficients. However, $A_t$ is not nilpotent
since
\begin{displaymath}
A_t^k=(1-t)^k E+(1-(1-t)^k)\mathbb{O}_n\neq \mathbb{O}_n
\end{displaymath}
for all $k\in\mathbb{N}$. 

For a positive integer, $k$, we denote by $M_{k}(\mathbb{R})$
the semigroup of all real matrices of size $k\times k$. Set
\begin{displaymath}
\mathbf{Q}_n=\{
(a_{i,j})\in M_n(\mathbb{R}):\sum_{k=1}^n a_{k,i}=
\sum_{k=1}^n a_{i,k}=1\text{ for all }i
\}.
\end{displaymath}
It is easy to check that $\mathbf{Q}_n$ is a subsemigroup 
of $M_n(\mathbb{R})$, moreover, that $\mathbb{O}_n$ is
the zero element of $\mathbf{Q}_n$.

\begin{lemma}\label{lemma2}
$A\in \mathbf{Q}_n$ if and only if the vector 
$\mathbf{v}=(1,1,\dots,1)$ is an eigenvector for $A$
with eigenvalue $1$, and the subspace
\begin{displaymath}
\mathbf{V}=\{(x_1,\dots,x_n)\in \mathbb{R}^n:
x_1+\dots+x_n=0\}
\end{displaymath}
is invariant with respect to $A$.
\end{lemma}

\begin{proof}
The ``only if'' part is checked by a direct calculation.
Let us prove the ``if'' part. The fact that $\mathbf{v}=(1,1,\dots,1)$ 
is an eigenvector for $A$ with eigenvalue $1$ means that
the sum of all elements in each row of $A$ equals $1$. Consider,
for $i\neq j$, the vector $u_{i,j}$ defined as follows:
the $i$-th and $j$-th coordinates of $u_{i,j}$ equal
$1$ and $-1$ respectively, and all other coordinates are
zero. Since $u_{i,j}\in \mathbf{V}$ we have 
$Au_{i,j}\in \mathbf{V}$, which means that in the matrix $A$
the sums of all elements in the $i$-th and in the $j$-th columns
coincide. From the first part we have the the sum of all 
entries of $A$ equals $n$. The lemma follows.
\end{proof}

\begin{proposition}\label{proposition1}
$\mathbf{Q}_n\cong M_{n-1}(\mathbb{R})$.
\end{proposition}

\begin{proof}
Let $\mathbf{v}$ and $\mathbf{V}$ be as in Lemma~\ref{lemma2}. 
Let further $\mathbf{v}_2,\dots,\mathbf{v}_n$ be an arbitrary
basis of the vector space $\mathbf{V}$, and $F$ be the transition
matrix from the standard basis $\mathbf{e}_1,\dots,\mathbf{e}_n$
of $\mathbb{R}^n$ to the basis $\mathbf{v},
\mathbf{v}_2,\dots,\mathbf{v}_n$. Then, from Lemma~\ref{lemma2}
it follows that for every $A\in \mathbf{Q}_n$ we have
\begin{displaymath}
F^{-1}AF=\left(\begin{array}{c|c} 1 & 0\\ \hline 0& B\end{array}
\right),
\end{displaymath}
where $B\in M_{n-1}(\mathbb{R})$. It is straightforward that
the map $A\mapsto B$ from $\mathbf{Q}_n$ to $M_{n-1}(\mathbb{R})$
is an isomorphism.
\end{proof}

In \cite[Section~7]{KM} it is shown that for an arbitrary field,
$\mathbb{F}$, there exists a bijection between the maximal nilpotent
subsemigroups of nilpotency class $k$ in the multiplicative
semigroup $M_{n}(\mathbb{F})$, and flags of length $k$ in
$\mathbb{F}^n$. Moreover, the maximal nilpotent subsemigroup
$T$ of nilpotency class $k$, which corresponds to the
flag $0=V_0\subsetneqq V_1\subsetneqq\dots\subsetneqq V_k=\mathbb{F}^n$,
has the following form:
\begin{displaymath}
T=\{A\in M_{n}(\mathbb{F}):A V_i\subset V_{i-1}\text{ for all }i\}.
\end{displaymath}
In particular, it follow that the nilpotency class of each maximal
nilpotent subsemigroup of $M_{n}(\mathbb{F})$ equals $n$, and that all
such subsemigroups are isomorphic and correspond bijectively to complete
flags in $\mathbb{F}^n$. Choose some basis $\mathbf{f}_1,\dots,
\mathbf{f}_n$ in $\mathbb{F}^n$ such that 
$V_i=\langle\mathbf{f}_1,\dots,\mathbf{f}_i\rangle$ for each $i$. Let
$F$ be the transition matrix from the standard basis
$\mathbf{e}_1,\dots,\mathbf{e}_n$ of $\mathbb{F}^n$ to the basis
$\mathbf{f}_1,\dots, \mathbf{f}_n$. Then the maximal nilpotent
subsemigroup $T$, which corresponds to the flag
$0=V_0\subsetneqq V_1\subsetneqq\dots\subsetneqq V_n=\mathbb{F}^n$,
has the form $F\mathrm{T}_nF^{-1}$, where $\mathrm{T}_n$ is the semigroup
of all upper triangular matrices from $M_{n}(\mathbb{F})$, having 
zero diagonal. This and Proposition~\ref{proposition1} immediately
implies:

\begin{theorem}\label{theorem3}
Let $\mathbf{v}$ and $\mathbf{V}$ be as in Lemma~\ref{lemma2}. Then:
\begin{enumerate}[(a)]
\item\label{theorem3.1} 
There exists a bijection between maximal nilpotent subsemigroups of
nilpotency class $k$ in the semigroup $\mathbf{Q}_n$ and flags of
length $k$ in the space $\mathbf{V}$.
\item\label{theorem3.2}
If $0=V_0\subsetneqq V_1\subsetneqq\dots\subsetneqq V_k=\mathbf{V}$
is a flag in $\mathbf{V}$, then the maximal nilpotent subsemigroup
of nilpotency class $k$, which corresponds to this flag, has the form
\begin{displaymath}
T=\{A\in \mathbf{Q}_n: AV_i\subset V_{i-1}\text{ for all }i\}.
\end{displaymath}
\item\label{theorem3.3}
The nilpotency class of each maximal nilpotent subsemigroup of
$\mathbf{Q}_n$ equals $n-1$, all such semigroups are isomorphic and
they correspond bijectively to complete flags in $\mathbf{V}$.
\item\label{theorem3.4}
Let $0=V_0\subsetneqq V_1\subsetneqq\dots\subsetneqq V_n=\mathbf{V}$
be a complete flag in $\mathbf{V}$ and $\mathbf{f}_2,\dots, \mathbf{f}_n$
be a basis in $\mathbf{V}$ such that 
$V_i=\langle\mathbf{f}_2,\dots,\mathbf{f}_{i+1}\rangle$ for each $i$.
Then the maximal nilpotent subsemigroup $T$, which corresponds to
this flag, has the form
\begin{displaymath}
T=F\left(\begin{array}{c|c}1 & 0\\ \hline 0& T_{n-1}\end{array}
\right)F^{-1},
\end{displaymath}
where $F$ is the transition matrix from the standard basis
$\mathbf{e}_1,\dots,\mathbf{e}_n$ of $\mathbb{R}^n$  to the 
basis $\mathbf{v},\mathbf{f}_2,\dots,\mathbf{f}_n$,
and $T_{n-1}$ is the semigroup of all upper triangular matrices
from $M_{n-1}(\mathbb{R})$ with zero diagonal.
\end{enumerate}
\end{theorem}

\begin{question}\label{qest}
Does the analogue of the semigroup $\mathbf{Q}_n$ for fields of 
positive characteristic $p$ have any
interesting properties in the case $p|n$?
\end{question}

\begin{lemma}\label{lemma3}
Let the semigroup $S_2$ be such that every nilpotent subsemigroup of
$S_2$ is contained in some maximal  nilpotent subsemigroup of $S_2$.
Let further $S_1$ be a subsemigroup of $S_2$. Then every maximal
nilpotent subsemigroup $T_1$ of $S_1$ of nilpotency class $k$
has the form $T_1=S_1\cap T_2$, where $T_2$ is some maximal
nilpotent subsemigroup of $S_2$ of nilpotency class $k$.
\end{lemma}

\begin{proof}
We have $T_1=S_1\cap T_2$, where $T_2$ is some maximal
nilpotent subsemigroup of $S_2$ of nilpotency class $k$, containing
$T_1$ (which exists because of our assumptions).
\end{proof}

\begin{lemma}\label{lemma4}
Let $T$ be a maximal nilpotent subsemigroup of nilpotency
class $k$ in $\mathbf{Q}_n$, $A\in \mathbf{Q}_n$, and $\alpha\neq 0$. 
Then $A\in T$ if and only if $\alpha A+(1-\alpha)\mathbb{O}_n\in T$.
\end{lemma}

\begin{proof}
Let $T'$ be a maximal nilpotent subsemigroup of nilpotency
class $k$ in $M_{n-1}(\mathbb{R})$. From the result of 
\cite[Section~7]{KM} mentioned above it follows that for any
$B\in M_{n-1}(\mathbb{R})$ and $\alpha\neq 0$ we have that
$B\in T'$ if and only if $\alpha B\in T'$. Since
\begin{displaymath}
\left(\begin{array}{c|c}1&0\\\hline 0&\alpha B\end{array}\right)=
\alpha\left(\begin{array}{c|c}1&0\\\hline 0&B\end{array}\right)+(1-\alpha)
\left(\begin{array}{c|c}1&0\\\hline 0&\alpha 0\end{array}\right),
\end{displaymath}
the necessary statement follows from Theorem~\ref{theorem3} and
the isomorphism $\mathbf{Q}_n\cong M_{n-1}(\mathbb{R})$ from the
proof of Proposition~\ref{proposition1}.
\end{proof}

\begin{theorem}\label{theorem4}
Let $k\in\mathbb{N}$. Then the map $\tau:T\mapsto T\cap \mathbf{D}_n$
defines a bijection between the maximal nilpotent subsemigroups of 
nilpotency class $k$ in $\mathbf{Q}_n$ and the maximal nilpotent 
subsemigroups of  nilpotency class $k$ in $\mathbf{D}_n$. In particular,
every nilpotent subsemigroup of $\mathbf{D}_n$ is contained in some
maximal nilpotent subsemigroup, and every maximal nilpotent subsemigroup
of $\mathbf{D}_n$ has nilpotency class $n-1$.
\end{theorem}

\begin{proof}
Taking Lemma~\ref{lemma3} into account, we have just to prove that 
$\tau$ is injective and preserves the nilpotency class. Let
$T_1$ and $T_2$ be two different maximal nilpotent subsemigroups
of nilpotency class $k$ from $\mathbf{Q}_n$ and let 
$A=(a_{i,j})\in T_2\setminus T_1$. Then
$\alpha A+(1-\alpha)\mathbb{O}_n\in T_2\setminus T_1$ for
any $\alpha\neq 0$ by Lemma~\ref{lemma4}. However, if
$\alpha\in\left(0,\min\{(2n|a_{i,j}|)^{-1}:a_{i,j}\neq 0\}\right)$, then we 
have $\alpha A+(1-\alpha)\mathbb{O}_n\in \mathbf{D}_n$. Hence
$(T_2\cap\mathbf{D}_n)\setminus (T_1\cap \mathbf{D}_n)\neq \varnothing$.
The injectivity of $\tau$ follows.

Let now $T$ be a maximal nilpotent subsemigroup of  nilpotency 
class $k$ in $\mathbf{Q}_n$, and $A_1,\dots,A_m\in T$ be such that
$A_1\cdots A_m\neq \mathbb{O}_n$. Then for any non-zero $\alpha_1,\dots,\alpha_m$
there exists $\beta\in \mathbb{R}$ such that
\begin{displaymath}
(\alpha_1 A_1+(1-\alpha_1)\mathbb{O}_n)\cdots
(\alpha_m A_m+(1-\alpha_m)\mathbb{O}_n)=
\alpha_1\cdots\alpha_m A_1\cdots A_m+\beta\mathbb{O}_n\neq 
\mathbb{O}_n.
\end{displaymath}
Analogously to the previous arguments, all $\alpha_i$ can be chosen such
that all corresponding $\alpha_i A_i+(1-\alpha_i)\mathbb{O}_n$ belong
to $\mathbf{D}_n$. Hence the nilpotency classes of $T$ and
$T\cap \mathbf{D}_n$ coincide. This completes the proof.
\end{proof}

\begin{corollary}\label{corollarynew1}
There is a bijection between the maximal nilpotent subsemigroups of
a given nilpotency class, $k$, in $\mathbf{D}_n$ and the flags of 
length $k$ in the $(n-1)$-dimensional real vector space $\mathbf{V}$.
In particular, the maximal nilpotent subsemigroups in $\mathbf{D}_n$
correspond to complete flags $0=V_0\subsetneqq V_1\subsetneqq\dots
\subsetneqq V_{n-1}=\mathbf{V}$.
\end{corollary}

\begin{proof}
This follows from Theorem~\ref{theorem3} and Theorem~\ref{theorem4}.
\end{proof}

The elements of the semigroup $T_{n-1}$ of all upper triangular matrices
with the zero diagonal can be naturally identified with the elements
of the vector space $\mathbb{R}^{n(n-1)/2}$. The condition that some
matrix from the nilpotent subsemigroup
$F\left(\begin{array}{c|c}1&0\\\hline 0&T_{n-1}\end{array}\right)F^{-1}$
belongs to $\mathbf{D}_n$ is equivalent to the condition of 
non-negativity of the elements of the matrix. If the transition matrix
$F$ is fixed, our condition reduces to a system of linear inequalities
for the coefficients of matrices in $T_{n-1}$. This means that
every maximal nilpotent subsemigroup $T$ of $\mathbf{D}_n$
corresponds to some convex polyhedron $\mathrm{P}(T)$ from
$\mathbb{R}^{n(n-1)/2}$.

\begin{proposition}\label{proposition2}
The polyhedron $\mathrm{P}(T)$ is bounded (compact) for every maximal
nilpotent subsemigroup $T$ of $\mathbf{D}_n$.
\end{proposition}

\begin{proof}
It is obvious that the point $O=(0,\dots,0)$, which corresponds to
the zero element $\mathbb{O}_n$ of the semigroup $T$ is an inner
point of $\mathrm{P}(T)$. Hence it is enough to show that the intersection
of every straight line from $\mathbb{R}^{n(n-1)/2}$, which contains
$O$, with $\mathrm{P}(T)$ is a bounded segment.

Let $A\in T\setminus\{\mathbb{O}_n\}$ and $M\in \mathrm{P}(T)$ be
the corresponding point of $\mathrm{P}(T)$. From the proof of
Lemma~\ref{lemma4} it follows that the elements of
the straight line $OM=\{\alpha M:\alpha\in\mathbb{R}\}$ correspond
to the elements of the subset
$\tilde{M}=\{\alpha A+(1-\alpha)\mathbb{O}_n:\alpha\in\mathbb{R}\}$
of $\mathbf{Q}_n$. Consider the intersection 
$\tilde{M}\cap \mathbf{D}_n$. Since $A\neq \mathbb{O}_n$, there exists
coefficients $a'$ and $a''$ of $A$ such that $a'>1/n$ and $a''<1/n$.
From the inequalities
\begin{displaymath}
\alpha a'+(1-\alpha)1/n\geq 0 \quad\text{ and }\quad
\alpha a''+(1-\alpha)1/n\geq 0
\end{displaymath}
we obtain $(1-na')^{-1}\leq \alpha\leq (1-na'')^{-1}$. Hence
$\tilde{M}\cap \mathbf{D}_n$ is a bounded segment and therefore
$OM\cap \mathrm{P}(T)$ is a bounded segment as well.
\end{proof}

We remark that the geometric structure of $\mathrm{P}(T)$
heavily depends on the choice of $T$:

\begin{example}\label{example1}
{\rm
Let $n=4$. We identify the matrix
$\left(\begin{array}{ccc}0&a&b\\0&0&c\\0&0&0\end{array}\right)$
from $T_3$ with the point $(a,b,c)\in \mathbb{R}^3$. Then for
the transition matrix
\begin{displaymath}
F'=\left(\begin{array}{cccc}1&1&1&1\\1&-1&0&0\\1&0&-1&0\\
1&0&0&-1\end{array}\right)
\end{displaymath}
the corresponding maximal nilpotent 
subsemigroup $T_{F'}$ of $\mathbf{D}_4$ is given by the following
system of linear inequalities:
\begin{displaymath}
\left\{
\begin{array}{l}
1+a+b+c\geq 0,\\
1-3a+b+c\geq 0,\\
1+a-3b-3c\geq 0,\\
1-a+3b\geq 0,\\
1+3a-b\geq 0,\\
1-a-b\geq 0,\\
1+3c\geq 0,\\
1-c\geq 0.
\end{array}
\right.
\end{displaymath}
The set of solutions to this system is the convex polyhedron
with vertexes $A=(-5/12,-1/4,-1/3)$, $B=(-1/4,-5/12,-1/3)$,
$C=(1/8,-7/24,-1/3)$, $D=(5/12,7/12,-1/3)$, $E=(1/4,3/4,-1/3)$,
$F=(-1/8,5/8,-1/3)$, $G=(-1/2,-1/2,0)$, $H=(1/2,1/2,0)$,
$I=(-1/2,-1/2,2/3)$, and $J=(1/2,-1/6,2/3)$. This polyhedron
has seven faces, one of which is a hexagon, two are
pentagons, two are quadrangles and the remaining two
are triangles (see Figure~\ref{figure1}).

\begin{figure}
\special{em:linewidth 0.4pt}
\unitlength 0.80mm
\linethickness{0.4pt}
\begin{center}
\begin{picture}(120.00,80.00)
\put(15.00,65.00){\makebox(0,0)[cc]{$\bullet$}}
\put(85.00,75.00){\makebox(0,0)[cc]{$\bullet$}}
\put(105.00,55.00){\makebox(0,0)[cc]{$\bullet$}}
\put(95.00,45.00){\makebox(0,0)[cc]{$\bullet$}}
\put(75.00,40.00){\makebox(0,0)[cc]{$\bullet$}}
\put(105.00,35.00){\makebox(0,0)[cc]{$\bullet$}}
\put(15.00,25.00){\makebox(0,0)[cc]{$\bullet$}}
\put(30.00,15.00){\makebox(0,0)[cc]{$\bullet$}}
\put(30.00,10.00){\makebox(0,0)[cc]{$\bullet$}}
\put(60.00,10.00){\makebox(0,0)[cc]{$\bullet$}}
\drawline(15.00,65.00)(85.00,75.00)
\drawline(15.00,65.00)(15.00,25.00)
\dashline{1}(15.00,65.00)(75.00,40.00)
\drawline(105.00,55.00)(85.00,75.00)
\drawline(60.00,10.00)(85.00,75.00)
\drawline(60.00,10.00)(105.00,35.00)
\drawline(60.00,10.00)(30.00,10.00)
\dashline{1}(30.00,15.00)(30.00,10.00)
\drawline(15.00,25.00)(30.00,10.00)
\dashline{1}(30.00,15.00)(15.00,25.00)
\drawline(105.00,55.00)(105.00,35.00)
\dashline{1}(95.00,45.00)(105.00,35.00)
\dashline{1}(105.00,55.00)(95.00,45.00)
\dashline{1}(75.00,40.00)(95.00,45.00)
\dashline{1}(75.00,40.00)(30.00,15.00)
\put(10.00,65.00){\makebox(0,0)[cc]{$I$}}
\put(10.00,25.00){\makebox(0,0)[cc]{$G$}}
\put(90.00,75.00){\makebox(0,0)[cc]{$J$}}
\put(110.00,55.00){\makebox(0,0)[cc]{$H$}}
\put(90.00,50.00){\makebox(0,0)[cc]{$E$}}
\put(110.00,35.00){\makebox(0,0)[cc]{$D$}}
\put(75.00,35.00){\makebox(0,0)[cc]{$F$}}
\put(30.00,20.00){\makebox(0,0)[cc]{$A$}}
\put(30.00,05.00){\makebox(0,0)[cc]{$B$}}
\put(60.00,05.00){\makebox(0,0)[cc]{$C$}}
\end{picture}
\end{center}
\caption{Polyhedron $ABCDEFGHIJ$.}\label{figure1}
\end{figure}
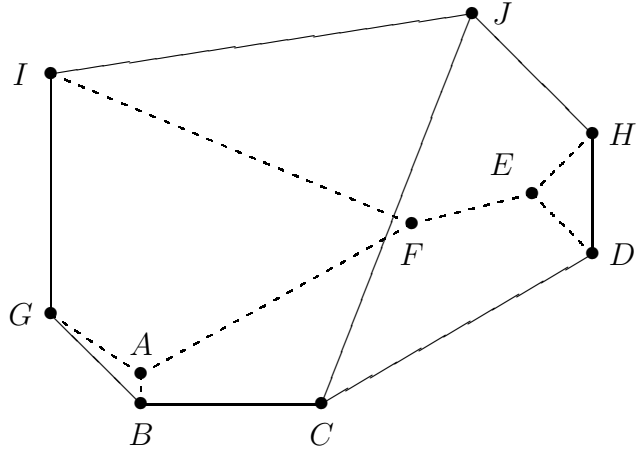

At the same time for the transition matrix
$F''=\left(\begin{array}{cccc}1&0&1&1\\1&0&1&0\\1&1&1&0\\
1&-1&-3&-1\end{array}\right)$ the corresponding maximal nilpotent 
subsemigroup $T_{F''}$ of $\mathbf{D}_4$ is given by the following
system of linear inequalities:
\begin{displaymath}
\left\{
\begin{array}{l}
1-a+4b+4c\geq 0,\\
1+3a-4b-4c\geq 0,\\
1+a-4b-12c\geq 0,\\
1-3a+4b+12c\geq 0,\\
1+4c\geq 0,\\
1-4c\geq 0,\\
1-a\geq 0,\\
1+a\geq 0.
\end{array}
\right.
\end{displaymath}
The set of solutions to this system is the tetrahedron
with the following vertexes: $X=(1,5/4,-1/4)$, $Y=(1,-1/4,1/4)$, 
$Z=(-1,-1/4,-1/4)$, and 
$W=(-1,-3/4,1/4)$.
}
\end{example}

\begin{question}
Is there any connection between the algebraic properties of 
a maximal nilpotent subsemigroup, $T\subset \mathbf{D}_n$,
and the geometric properties of the polyhedron
$\mathrm{P}(T)$?
\end{question}

\vspace{1cm}

\noindent
O.G.: Department of Mechanics and Mathematics, Kyiv Taras Shevchenko
University, 64, Volodymyrska st., 01033, Kyiv, UKRAINE,\\
e-mail: {\tt ganiyshk\symbol{64}univ.kiev.ua}
\vspace{0.5cm}

\noindent
V.M.: Department of Mathematics, Uppsala University, Box 480,
SE 751 06, Uppsala, SWEDEN, e-mail: {\tt mazor\symbol{64}math.uu.se},\\
web: ``http://www.math.uu.se/$\tilde{\hspace{2mm}}$mazor''
\vspace{0.5cm}

\end{document}